\documentclass{amsart}

\usepackage{amscd}
\usepackage{enumitem}
\usepackage{graphicx}

\newtheorem{theorem}{Theorem}[section]

\theoremstyle{definition}
\newtheorem{definition}[theorem]{Definition}
\newtheorem{example}[theorem]{Example}

\newtheorem{corollary}[theorem]{Corollary}
\newtheorem{proposition}[theorem]{Proposition}
\theoremstyle{remark}
\newtheorem{remark}[theorem]{Remark}
\newtheorem{conjecture}[theorem]{Conjecture}
\newtheorem{problem}[theorem]{Problem}

\numberwithin{equation}{section}

\begin{document}

\title{ On 
 coverings by Minkowski balls in the plane and a duality}

\author{Nikolaj M. Glazunov}
\address{Institute of Mathematics and Informatics Bulgarian Academy of Sciences, 1113 Sofia, Bulgaria.
Glushkov Institute of Cybernetics NASU, Kiev, Ukraine}
\curraddr{Institute of Mathematics and Informatics Bulgarian Academy of Sciences, 1113 Sofia, Bulgaria}
\email{glanm@yahoo.com}

\thanks{The author was supported   Simons grant 992227}





\begin{abstract}

Lattice coverings in the real plane by Minkowski balls  are studied. We exploit the duality of admissible lattices of Minkowski balls and
 inscribed convex symmetric hexagons of these  balls. An explicit moduli space  of the areas of these hexagons is constructed, giving the values of the determinants of corresponding  covering lattices. Low and upper bounds for covering constants of Minkowski balls are given. The best known value of covering density of Minkowski ball is obtained. One conjecture   and  one problem are formulated.

\end{abstract}

\maketitle

\section{Introduction}
A system of equal  balls in $n$-dimensional space is said to form a covering,  if each point in space belongs to at least one ball of this system.

  Let
\begin{equation}
\label{smt}
   D_p: \;  |x|^p + |y|^p \le 1, \; p \ge 1
\end{equation}
be (two-dimensional) Minkowski's balls with boundary 
\begin{equation}
\label{mt}
C_p: \; |x|^p + |y|^p = 1, \;  p \ge 1.
\end{equation}

Lattice coverings in the real plane by Minkowski balls  with inscribed convex symmetric hexagons generated by admissible lattices of Minkowski balls are studied. The covering problem is dual to the packing problem. In different covering problems, this duality, apparently, can manifest itself in different ways. In this paper, duality is understood as the duality between admissible lattices having 3 pairs of points on the boundary $D_p$, and symmetric convex hexagons (possibly a quadrangle) corresponding to these lattices and determining the values of the determinants of the corresponding covering lattices.

Geometrically, the set of Minkowski's balls is a solid tube of variable diameter with square 
sections at the ends $p = 1$, and $p = \infty$ which we call the limiting Minkowski's balls.
We use admissible lattices that have 3 pairs of points on the boundary of $D_p$.
We construct an explicit moduli space (\ref{Ap})  of areas of symmetric convex hexagons (possibly a quadrangle) which determine values of the determinants of corresponding  covering lattices .

This construction and its consequences make it possible to obtain the best currently known lower bounds for covering constants of Minkowski balls. An estimate from above is also given.
These give the corresponding estimation of the covering densities.
The best known value of covering density of Minkowski ball is obtained (Proposition \ref{emcd}). 

Our investigations connect with Minkowski conjecture~\cite{Mi:DA,M:LP,D:NC,Co:MC,W:MC,GGM:PM}  and use 
results of its investigation. This conjecture was formulated by Minkowski while studying a problem
 of Diophantine approximations for sums of powers of linear forms~\cite{Mi:DA}. Mordell tested and confirmed this conjecture for $p=4$~\cite{M:LP}.
 Davis adjusted Minkowski conjecture~\cite{D:NC} for small values of $1<p \le p_0$ (see theorems (\ref{ggm}), (\ref{ma})).
 Minkowski (implicitly) and Cohn~\cite{Co:MC} explicitly defined an appropriate moduli space for computing determinants of critical lattices.
 We follow these researchers, albeit in a slightly different sequence.
 Our research on the packing problem~\cite{Gl5,Gllpop} was also used.  Corresponding results and conjectures are stated in the simpliest way in terms of geometric lattices and covering lattices~\cite{Mi:DA,Cassels,lek}. 

Based on these studies, one conjecture  about the behavior of the covering constants of Minkowski balls (see \ref{bcc})  and  one problem  on belonging  of the curve of maxima of covering constant  to a family introduced in the work (see \ref{mbf}) are formulated.

 \section{Two-dimensional Minkowski balls, admissible lattices and hexagons}
\label{Minkowski balls}

\begin{proposition}
 The volume of  Minkowski ball $D_p$ is equal $4 \frac{(\Gamma(1 + \frac{1}{p}))^2}{\Gamma(1 + \frac{2}{p})}$.
\end{proposition}
{\bf Proof}. (by Minkowski)
 Let $x^P + y^p \le 1, x \ge 0, y \ge 0.$
Put $x^p = \xi, y^p = \eta.$
\begin{equation}
 \label{vdp}
V(D_p) = \frac{4}{p^2}\iint \xi^{\frac{1}{p} - 1}\eta^{\frac{1}{p} - 1} d\xi d\eta,
\end{equation}
where the integral extends to the area
\[
  \xi + \eta \le 1, \; \xi \ge 0, \; \eta \ge 0.
\]
 Expression (\ref{vdp}) can be represented in terms of Gamma functions, and we get
\[
 V(D_p) = 4 \frac{(\Gamma(1 + \frac{1}{p}))^{2}}{\Gamma(1 + \frac{2}{p})}.
\]

\subsection{Admissible lattices}~\cite{Cassels,lek}.

Let $\mathcal R $ be a set and $\Lambda $ be a  lattice with base 
\[
\{a_1, \ldots ,a_n \}
\] 
in ${\mathbb R}^n.$ \\

A lattice $\Lambda $
is {\it admissible} for body $\mathcal R $ 
($ {\mathcal R}-${\it admissible})
if ${\mathcal R} \bigcap \Lambda = \emptyset $ or $0.$\\
 
 Let $\overline {\mathcal R}$ be the closure of $ {\mathcal R}$.
A lattice $\Lambda $
is {\it strictly admissible} for  $\overline {\mathcal R}$  
( $\overline {\mathcal R}-${\it strictly  admissible})
if $\overline {\mathcal R} \bigcap \Lambda = \emptyset $ or $0.$\\

Let
 \begin{equation}
 \label{det} 
  d(\Lambda) = |\det (a_1, \ldots ,a_n )|
\end{equation}
   be the determinant of 
$\Lambda.$ \\

\subsection{Hexagons and coverings}

Every admissible lattice of $D_p$ containing 3 pairs of points on the boundary of $D_p$ defines a hexagon inscribed in $D_p$.

We denote such a hexagon as $ {\mathcal{H}} _p $, and call it a hexagon of the admissible lattice, briefly: al-hexagon.

 \begin{remark}
In the limiting cases $p=1$ and $p=\infty$ corresponding hexagons (following  Fejer Toth, see \ref{sth} ) $ {\mathcal{H}} _p $ are quadrangles (see Proposition \ref{lmt}).
 \end{remark}
 
 By specifying results from \cite{lek} to the case of Minkowski balles and  al-hexagons $ {\mathcal{H}} _p $ we have the next
 \begin{proposition}
 \label{ccdmb}
 Each Minkowski ball $D_p$ contains an al-hexagon  (possibly a quadrangle) $ {\mathcal{H}} _p $ of maximum area.
 If this area is denoted by $\gamma_h(D_p)$ (or for symmetric convex inscribed al-hexagons  
 $ {\mathcal{H}} _p $ by $ a( {\mathcal{H}} _p)$)  then 
  \begin{equation}
 \label{ccmb} 
       \Gamma(D_p) = \gamma_h(D_p).
       \end{equation}  
       For the density of thinnest lattice covering of $D_p$ we get
       \begin{equation}
 \label{dcmb} 
       \vartheta(D_p) = V(D_p)/\gamma_h(D_p).
 \end{equation}  
\end{proposition} 

 \begin{remark}
The expression (\ref{ccmb}) is the ${\bf covering \; constant}$ of  Minkowski ball $D_p$, and 
 the expression (\ref{dcmb}) is the ${\bf density}$ of the  ${\bf thinnest \; lattice \;covering}$ within the framework of our considerations.
   \end{remark}
   
   \begin{example} 
\label{ex}

1. The lattice 
\[
  \Lambda_{1}^{(1)} =  \{1, 1), (2, 0)\}
\]
is the covering lattice by $D_1$. 
The limiting case of Minkowski balls with $p = 1$ gives an optimal covering (indeed a partition) of density $1$.

2.   The lattice 
\[
  \Lambda_{\infty}^{(1)} =  \{(2, 2), (0, 2)\}
 \]
is the  covering lattice by $D_{\infty}$. 
Limiting case of Minkowski balls  for $ p = \infty$ gives  the optimal covering (indeed a partition) of density $1$.

3. The lattice 
\[
\Lambda_{2}^{(0)} =  \{(\sqrt{3}, 0), (\frac{\sqrt{3}}{2}, \frac{3}{2})\}
\]
is the  covering lattice by $D_2$. 
The lattice gives the  optimal covering of density $\vartheta = \frac{2\pi}{3\sqrt{3}}  \approx 1,19$.
\end{example}

\section{Moduli spaces,  moduli problems and   and areas of inscribed hexagons} 

 Moduli spaces are spaces of solutions of geometric classification problems and  dynamics of varying objects in families \cite{KM,hm}.\\

Example: moduli space for Legendre elliptic curve: $$y^2 = x(x-1)(x- \lambda), \; \lambda \in {\mathbb C}$$\\

Moduli problem - try to obtain solution for all objects of the moduli space  or for objects of some subspace of the corresponding moduli space\\

Currently, there is great interest in unifying problems of algebraic, arithmetic and geometric classification and dynamics through moduli spaces of object classes and studying related problems within such unifications.
Various types of moduli spaces were presented and studied in the famous works of Eichler, Lazard, Shimura, Siegel, Igusa, Deligne-Mumford-Knudsen, Deligne-Rapoport, Katz-Mazur, Faltings, Drinfeld and others.
But the author is not aware of applications of moduli spaces to covering problems.\\

 In present our research the moduli space is real differentiable moduli space with exceptional points. \\ 

 \subsection{ Solid Minkowski tube,  Minkowski tube and  metrics}
 \label{sth}
 \begin{definition}
  We will call (\ref{smt}) a solid Minkowski tube and (\ref{mt}) a Minkowski tube.
 \end{definition}
 
  \begin{remark}
Expressions
\begin{equation}
\label{met}
    (|x|^p + |y|^p)^{1/p}
    \end{equation}
    are metrics. They also form the tube of metrics for $1 \le p \le \infty$. But we won't use it.
\end{remark}
 
 \begin{proposition}
 \label{lmt}
A Minkowski tube (\ref{mt}) is a real manifold which, for $p=1$ and $p=\infty$, has exceptional points corresponding to the corners of the corresponding limit quadrangles.
 \end{proposition}
 {\bf Proof}.  Follows from \cite{Mi:DA} and from the definition of Minkowski tube.\\
 
The use of moduli spaces and the study of moduli problems require the adjustment of some well-known concepts.\\


Following  Fejer Toth, a closed finite convex region bounded by at most six straight line segments will be called a hexagon.

This definition of a hexagon was mentioned by Bambah, Rogers and other famous covering researchers \cite{Ro:PC}.\\

\subsection{Inscribed hexagons of admissible lattices  of Minkowski balls $D_p$ and their moduli space }
 
 \begin{theorem}
 \label{msa}
The set of areas of the entire family of al-hexagons (possibly a quadrangle) $ {\mathcal{H}}_p, (1 <  p < \infty)$ inscribed in Minkowski balls $D_p$ is parameterized by the function $A(\sigma, p)$ defining the corresponding moduli space
   $\mathbf{A}: $
 \begin{equation}
 \label{Ap}
 A(\sigma, p) = 3(\tau + \sigma)(1 + \tau^{p})^{-\frac{1}{p}}
  (1 + \sigma^p)^{-\frac{1}{p}}, 
 \end{equation}
 In limiting cases we have the following areas 
 of quadrangles: $ p=1, \; a( {\mathcal{H}} _1) = 2; \; p=\infty, \; a( {\mathcal{H}} _{\infty}) = 4$.
  \end{theorem}
    {\bf Proof.}  
  The area $a( {\mathcal{H}} _p)$ of an  al-hexagon is given (see (\ref{det})) by 
   \begin{equation}
   a( {\mathcal{H}} _p) = 2\cdot d\left( \begin{array}{cc} x_1 & y_1 \\
x_1 - x_2 & y_1 + y_2 \end{array} \right) + d\left( \begin{array}{cc} -x_2 & y_2 \\
- x_1 & -y_1  \end{array} \right) =
  3(x_1 y_2 + x_2 y_1),
\end{equation}  
 where $x_1 > x_2 \ge 0, y_1 \ge 0,  y_2 > 0 $  and 
  $x_1^p + y_1^p =  |-x_2|^p + y_2^p =  (x_1-x_2)^p + (y_1 + y_2)^p = 1, \;  1 < p < \infty $.
   
 The parameterization of admissible lattices (al) is implicitly given in \cite{Mi:DA} and explicitly specified in \cite{Co:MC,W:MC,GGM:PM,Gl4,Gl5}.We specify the parameterization   of al-hexagons.
  Let 
 \[
 0 \le  \tau < \sigma , \; 0 \le \tau \le \tau_p.
  \]
   $\tau_{p}$ is defined by the
equation $ 2(1 - \tau_{p})^{p} = 1 + \tau_{p}^{p}, \; 0 \leq
\tau_{p} < 1. $ 
\[
  1 \le\sigma \le \sigma_p, \; \sigma_p = (2^p - 1)^{\frac{1}{p}}.
\]
  For inscribed in $D_p$ a hexagon $ {\mathcal{H}} _p$ with points $(x_1, y_1) = ((1 + \tau^{p})^{-\frac{1}{p}}, \tau(1 + \tau^{p})^{-\frac{1}{p}})$, $(-x_2, y_2) = (-(1 + \sigma^p)^{-\frac{1}{p}}, \sigma(1 + \sigma^p)^{-\frac{1}{p}})$,
  $(x_1 - x_2,  y_1 + y_2 ) = (1 + \tau^{p})^{-\frac{1}{p}} - (1 + \sigma^p)^{-\frac{1}{p}},  \tau(1 + \tau^{p})^{-\frac{1}{p}} + \sigma(1 + \sigma^p)^{-\frac{1}{p}})$.\\
  This gives the formula (\ref{Ap}).
  Limiting cases follow from the proof of Minkowski conjecture (Theorem \ref{ggm}) (see also 1 and 2 of  Example \ref{ex}).
  
  \begin{definition}
    Moduli space (\ref{Ap}) is the set of areas of the entire family of al-hexagons (possibly a quadrilateral) $ {\mathcal{H}}_p, (1 \le p \le \infty)$ inscribed in Minkowski balls $D_p$. Briefly we will call (\ref{Ap}) the moduli space of hexagons.
  \end{definition}
  
   \begin{remark}
    The function domain of the  function (\ref{Ap}) is
     $$ {\mathcal M}: \; \infty > p > 1, \; 1 \leq \sigma \leq \sigma_{p} =
 (2^p - 1)^{\frac{1}{p}} $$
     $ {\mathcal M}$  has two limit points $p=1$ and $p=\infty$.
       \end{remark}
  
  \begin{remark}
  The function (\ref{Ap}) is upwardly convex (concave)  in the  real space with Cartesien coordinates     $(\sigma, p, A)$.
  \end{remark}

   \section{Sections of  moduli space of hexagons}
   
     We will solve the problem of the thinnest covering of the real plane by Minkowski balls with hexagons $ {\mathcal{H}} _p$ if we define a curve of maxima on an upwardly convex (concave) surface (\ref{Ap}).
Unfortunately, to date the author knows only a few points of maxima on this curve.   In this regard, we present here a family of surfaces, give  some intersections of (\ref{Ap}) with these surfaces and formulate a (rough) hypothesis about the curve of maxima on (\ref{Ap}).

   \begin{definition}
     Let 
      \begin{equation}
 \label{fsap}  
   \sigma_{\alpha, p} = (2^p - 1)^\frac{1}{\alpha p}, \;  \alpha \ge 1
      \end{equation}  
     be a family of curves parameterized by ${\alpha }$ with $p$ varying from $1$ to $\infty$.     
       in 
        real plane  $(\sigma, p)$.
   \end{definition}
   
     \begin{remark}
     $\sigma_{\alpha, p} \neq  \sigma_{ p, \alpha}.$ 
    \end{remark}
   
  \begin{remark}
     Curves of the family (\ref{fsap})
     define family of curves
          \begin{equation}
 \label{sap}  
     A(\alpha, p): 3(\tau + (2^p - 1)^\frac{1}{\alpha p})(1 + \tau^{p})^{-\frac{1}{p}}
  (1 + ((2^p - 1)^\frac{1}{\alpha p})^p)^{-\frac{1}{p}} 
  \end{equation}  
  on the function range (\ref{Ap}).
   These curves are parametrized by $\alpha$. For each such $\alpha, \; p  $ varies from $ 1$ to   $\infty \;(1 < p < \infty)$.
   \end{remark}

  \section{ Estimates and Bounds on Covering Constants }
  
 Here we estimate covering constants of Minkowski balls from above and from below. 
According to the results of Sas and Dovker (see \cite{lek}) in our case of Minkowski balls we have

\begin{proposition}
\label{mhe}
Let  ${\mathcal H}_p$ be the 
 hexagon of maximuml area inscribed in Minkowski ball $D_p$.\\
 Then $V({\mathcal H}_p) = \gamma_h(D_p) \ge \frac{3\sqrt{3}}{2\pi}\cdot V(D_p)$
\end{proposition}

\begin{remark}
\label{teamb}
The trivial estimate from above of $V({\mathcal H}_p) $ is given by the area $V(D_p)$ of Minkowski ball:
$$ V({\mathcal H}_p) \le  V(D_p) $$.
\end{remark}

Let estimate covering constants of Minkowski balls $D_p$ from below on the base of results of the proof of Minkowski conjecture.
To do this, we will use the results of the proof of the Minkowski conjecture (now
theorem \cite{GGM:PM}) and the constructed moduli space (\ref{Ap}).
In this regard recall at first some definitions and  results.\\

The infimum
$\Delta(\mathcal R) $ of determinants of all lattices admissible for
$\mathcal R $ is called {\em the critical determinant} 
of $\mathcal R; $
if there is no $\mathcal R-$admissible lattices then puts
$\Delta(\mathcal R) = \infty. $ \\

 A lattice 
$\Lambda $ is {\em critical}
if $ d(\Lambda) = \Delta(\mathcal R).$ \\

\begin{theorem}
\cite{GGM:PM}
\label{ggm}
$$\Delta(D_p) = \left\{
                   \begin{array}{lc}
    \Delta(p,1), \; 1 \le p \le 2, \; p \ge p_{0},\\
    \Delta(p,\sigma_p), \;  2 \le p \le p_{0};\\
                     \end{array}
                       \right.
                           $$
here $p_{0}$ is a real number that is defined unique by conditions
$\Delta(p_{0},\sigma_p) = \Delta(p_{0},1),  \;
2,57 < p_{0}  < 2,58, \; p_0  \approx 2.5725 $
\end{theorem}
\begin{remark}
We will call $p_{0}$ the Davis constant.
\end{remark}

\begin{corollary}

From theorem (\ref{ggm}) in notations \cite{GGM:PM,Gl4} we have next expressions for critical determinants and their lattices:  
\begin{enumerate} 

\item  ${\Delta^{(0)}_p} = \Delta(p, {\sigma_p}) =  \frac{1}{2}{\sigma}_{p},$ 

\item $ {\sigma}_{p} = (2^p - 1)^{1/p},$

\item  ${\Delta^{(1)}_p}  = \Delta(p,1) = 4^{-\frac{1}{p}}\frac{1 +\tau_p }{1 - \tau_p}$,   

\item  $2(1 - \tau_p)^p = 1 + \tau_p^p,  \;  0 \le \tau_p < 1.$  

\end{enumerate}
\end{corollary}
  For their critical lattices respectively  $\Lambda_{p}^{(0)},\; \Lambda_{p}^{(1)}$ next conditions satisfy:   $\Lambda_{p}^{(0)}$ and 
 $\Lambda_{p}^{(1)}$  are  two $D_p$-admissible lattices each of which contains
three pairs of points on the boundary of $D_p$  with the
property that 
\begin{itemize}

\item $(1,0) \in \Lambda_{p}^{(0)},$

\item $(-2^{-1/p},2^{-1/p}) \in \Lambda_{p}^{(1)},$

\end{itemize}
 (under these conditions the lattices are
uniquely defined).\\

We specify the results of the Theorem \ref{ggm} for the moduli space (\ref{Ap}).

\begin{theorem}
\label{ma}
\begin{equation}
\label{ema}
\min({\bf A}({\mathcal H}_p)) = \left\{
                   \begin{array}{lc}
    3\cdot4^{-1/p} \frac{1 +\tau_p }{1 - \tau_p}, \; 1 \le p \le 2, \; p \ge p_{0},\\
    \frac{3}{2}{\sigma}_{p}, \;  2 \le p \le p_{0};\\
                     \end{array}
                       \right.
                       \end{equation}
\end{theorem}

To estimate the covering constants of Minkowski balls $D_p$ from below, we define the inverse minimum 
$i-\min({\bf A}({\mathcal H}_p)) $ for  ${\bf A}({\mathcal H}_p) $. 

\begin{definition}
\begin{equation}
 \label{ima}
    i-\min({\bf A}) =   i-\min({\bf A}({\mathcal H}_p)) = \left\{
                   \begin{array}{lc}
              \frac{3}{2}{\sigma}_{p}, \;   1 \le p \le 2, \; p \ge p_{0},  ;\\      
    3\cdot4^{-1/p} \frac{1 +\tau_p }{1 - \tau_p}, \; 2 \le p \le p_{0}.\\
                        \end{array}
                       \right.
                       \end{equation}
\end{definition}

\begin{example}
    In the case $p=3$ for moduli space $\mathbf{A}$ (see \ref{Ap}) we have \\
    $ \min({\bf A}({\mathcal H}_3))     \approx      2.859     $ \\
     $ i-\min({\bf A})({\mathcal H}_3))        \approx     2.870       $\\
     $  \min({\bf A}({\mathcal H}_3)) <  i-\min({\bf A})({\mathcal H}_3)). $
 \end{example}

From theorems \ref{msa}, \ref{ma} with definition  \ref{ma} we have the next estimation from below of the covering constant $\Gamma(D_p)$.

\begin{proposition}
\begin{equation}
 \label{eccb}
 \Gamma(D_p)  \ge i-\min({\bf A}).
\end{equation}
\end{proposition}

\begin{remark}
  For $p > 2$ the estimate (\ref{eccb}) is better than the estimate of the Proposition \ref{mhe}.
\end{remark}

\section{ Estimates and calculation  of Covering Density }

At first give the estimation of the thinnest lattice covering of $D_p$ from above on the base of  results by Sas and Dovker (\cite{lek}).

\begin{proposition}
 The density  $\vartheta$ of any covering of  ${\mathbb R}^2$ by  Minkowski balls $D_p$ satisfies
  \begin{equation}
 \label{deu}
    \vartheta \le \frac{2\pi}{3\sqrt{3}}  \approx 1,19;
\end{equation}
    so 
     \begin{equation}
 \label{deu3}
    \vartheta(D_3) \le   1,19;
\end{equation}
 \end{proposition}
  
  {\bf Proof.} By Proposition \ref{ccdmb} $\Gamma(D_p) = \gamma_h(D_p)$; hence follows ( \ref{deu}).\\
  
  Let us now examine the curves of the family (\ref{sap}) from the point of view of the coverage densities they determine.
Here we do this only for one point on one curve.
 But even such example gives earlier unknown minimum density for Minkowski ball.

Take the curve from the family (\ref{sap} ) defined by the function
\begin{equation}
 \label{s2p}
  A(2,p) = 3(\tau + (2^p - 1)^\frac{1}{2 p})(1 + \tau^{p})^{-\frac{1}{p}}
  (1 + ((2^p - 1)^\frac{1}{2 p})^p)^{-\frac{1}{p}}. 
 \end{equation} 
 
 \begin{proposition}
 \label{emcd}
  The convex symmetric hexagon defined by the point $p =3$ on the curve (\ref{s2p}) gives a covering density $\vartheta(D_3) \approx 1.0567$ by the Minkowski ball $D_3$, the minimum for known densities of Minkowski balls at $1 < p < \infty$.
 \end{proposition}
 {\bf Proof}. \\
 ${\sigma}_{2,3} = (2^{3}-1)^\frac{1}{2\cdot3} \approx 1.383;$\\
 $\tau_3 \approx 0.20406, 0 \le \tau \le \tau_3,$ for curve (\ref{s2p}) take $\tau = 0.12;$\\
 $V(D_3) \approx 3.52, \; \gamma_h(D_3) \approx 3.331;  $ \\
 So $\vartheta(D_3) = V(D_3)/\gamma_h(D_3)   \approx 1.0567;$ \\
 
      We state here one conjecture and one problem which arise naturally from our work.
 
 \begin{conjecture}     
 \label{bcc}
 The curve of maxima of covering constants increases from $p=1$ to $p=2$ and decreases from $p=2$ to $p=\infty$.     
\end{conjecture}    
  
\begin{problem}
\label{mbf}
Does the curve of maxima of covering constants belong to the family (\ref{sap})?
\end{problem}

\thanks{{\bf Acknowledgments} 
The author would like to warmly thank Prof. P. Boyvalenkov for his support}.

\bibliographystyle{amsplain}

\end{document}